\documentclass[12pt]{amsart}

\usepackage{color}
\definecolor{darkgreen}{rgb}{0, 0.40, 0}

\usepackage{amsthm} 
\usepackage{amsmath} 
\usepackage{amssymb} 

\usepackage{psfrag} 
\usepackage{epsfig} 

\usepackage[
pdfborderstyle={},     
ps2pdf,
colorlinks=true,
linkcolor=blue,
citecolor=darkgreen,
urlcolor=magenta
]{hyperref}

\newcommand{\calC}{\mathcal{C}}

\newcommand{\calE}{\mathcal{E}}

\newcommand{\calJ}{\mathcal{J}}

\newcommand{\calT}{\mathcal{T}}
\newcommand{\calX}{\mathcal{X}}
\newcommand{\calY}{\mathcal{Y}}

\renewcommand{\setminus}{{\smallsetminus}}

\newcommand{\st}{\mathbin{\mid}} 
\newcommand{\from}{\colon} 

\newcommand{\isom}{{\mathrel{\cong}}} 

\newcommand{\bdy}{\partial} 

\newcommand{\MCG}{\mathcal{MCG}} 

\newcommand{\Aut}{\operatorname{Aut}} 
\newcommand{\Out}{\operatorname{Out}} 
\newcommand{\Inn}{\operatorname{Inn}} 

\newcommand{\Teich}{{Teichm\"uller~}}

\newcommand{\refsec}[1]{Section~\ref{Sec:#1}}
\newcommand{\refthm}[1]{Theorem~\ref{Thm:#1}}
\newcommand{\refcor}[1]{Corollary~\ref{Cor:#1}}
\newcommand{\reflem}[1]{Lemma~\ref{Lem:#1}}

\newcommand{\refrem}[1]{Remark~\ref{Rem:#1}}

\newcommand{\refprob}[1]{Problem~\ref{Prob:#1}}

\theoremstyle{plain}
\numberwithin{equation}{section} 
\newtheorem{theorem}[equation]{Theorem}
\newtheorem{corollary}[equation]{Corollary}
\newtheorem{lemma}[equation]{Lemma}

\newtheorem{proposition}[equation]{Proposition}

\theoremstyle{definition}
\newtheorem{definition}[equation]{Definition}
\newtheorem{remark}[equation]{Remark}
\newtheorem*{remark*}{Remark}

\newtheorem*{claim*}{Claim}
\newtheorem{problem}[equation]{Problem}

\newtheorem*{question*}{Question}
\newtheorem*{answer*}{Answer}
\newtheorem*{application*}{Application}

\newcommand{\fakeenv}{} 

\newenvironment{restate}[2]  
{ 
 \renewcommand{\fakeenv}{#2} 
 \theoremstyle{plain} 
 \newtheorem*{\fakeenv}{#1~\ref{#2}} 
 \begin{\fakeenv}
}
{
 \end{\fakeenv}
}

\newcommand{\Sig}{\Sigma}
\newcommand{\ep}{\epsilon}

\newcommand{\Thresh}{{\sf T}}
\newcommand{\Acc}{{\sf A}}
\newcommand{\Bound}{{\sf M}}
\newcommand{\Error}{{\sf C}}
\renewcommand{\ep}{{\sf e}}

\newcommand{\Progress}{{\sf P}}
\newcommand{\BackTrack}{{\sf B}}
\newcommand{\Length}{{\sf L}}
\newcommand{\Comp}{{\sf K}}
\newcommand{\Degree}{{\sf d}}
\newcommand{\Const}{{\sf c}}
\newcommand{\Quasi}{{\sf Q}}
\newcommand{\NoProgress}{{\sf N}}

\newcommand{\LongThick}{\Gamma(\ep, \Length)}

\newcommand{\punS}{S^{\circ}}
\newcommand{\punSig}{\Sigma^{\circ}}
\newcommand{\punSigma}{\Sigma^{\circ}}
\newcommand{\punPsi}{\Psi^{\circ}}

\newcommand{\orb}{{\rm orb}}

\begin{document}

\title{Covers and the curve complex}

\author{Kasra Rafi}
\email{rafi@math.uconn.edu}
\urladdr{http://www.math.uconn.edu/$\sim$rafi}

\author{Saul Schleimer}
\email{saulsch@math.rutgers.edu}
\urladdr{http://www.math.rutgers.edu/$\sim$saulsch}

\thanks{This work is in the public domain.}

\date{\today}

\begin{abstract}
We provide the first non-trivial examples of quasi-isometric
embeddings between curve complexes.  These are induced either by
puncturing a closed surface or via orbifold coverings.  As a
corollary, we give new quasi-isometric embeddings between mapping
class groups.
\end{abstract}

\maketitle

\section{Introduction}

The complex of curves~\cite{Harvey81} arises in the study of
mapping class groups, \Teich spaces, and three-manifolds. 
We follow Masur and Minsky \cite{MasurMinsky99, MasurMinsky00} in
studying the coarse geometry of the curve complex.

Since the curve complex is locally infinite many of the standard
quasi-isometry invariants
are of questionable utility.
Instead, we concentrate on a different family of invariants:
metrically natural subspaces.  Note that the well-known subspaces of
the curve complex, such as the complex of separating curves, the disk
complex of a handlebody and so on, are {\em not} quasi-isometrically
embedded and so do not give invariants in any obvious way.  We
therefore restrict our attention to:

\begin{problem}
\label{Prob:Classify}
Classify all quasi-isometric embeddings between curve complexes.
\end{problem}

A special case of this question, answered in our forthcoming
paper~\cite{RafiSchleimer07}, is the computation the quasi-isometry
group of the curve complex.  More generally, one may ask for a theory
of quasi-isometric embeddings between combinatorial moduli spaces.

\begin{problem}
\label{Prob:Classify2}
Classify all quasi-isometric embeddings between mapping class groups
(respectively pants complexes, Hatcher-Thurston complexes, and so on).
\end{problem}

Again, an important special case of \refprob{Classify2} is the
computation of the quasi-isometry group of the mapping class group;
this has recently been claimed by Behrstock, Kleiner, Minsky, and
Mosher as well as by Hamenst\"adt~\cite{Hamenstaedt05b}.



This paper discusses the first non-trivial examples of one curve
complex being quasi-isometrically embedded in another.  These arise
geometrically, either by puncturing a closed surface
(\refthm{Puncture}) or from orbifold covering maps (\refthm{Cover}).




The puncturing construction, discussed in \refsec{Statements}, is
inspired by Harer's paper~\cite[Lemma~3.6]{Harer86}.  This
construction is straight-forward but, surprisingly, gives an
uncountable number of {\em isometric} embeddings of between curve
complexes.  It seems unlikely that the puncturing map gives rise to a
quasi-isometric embedding of mapping class groups.

The second construction, based on coverings, is more difficult and
relies crucially on the fact that a covering induces an isometric
embedding of \Teich spaces.
As an application, we prove \refthm{MCG}: any orbifold covering map
induces a quasi-isometric embedding of the associated mapping class
groups.  We then have \refcor{Luis}: there is a quasi-isometric
embedding of the spherical braid group on $2g+2$ strands into the
mapping class group of genus $g$; this answers a question of Luis
Paris.

\subsection*{Acknowledgments} We thank Jason Behrstock, Jason Manning,
Dan Margalit and Mahan Mj.~for their comments on an early version of
this paper.

\section{Statements}
\label{Sec:Statements}

Suppose that $\Sigma$ is a compact orientable orbifold of dimension
two.  For definitions and discussion of orbifolds we refer the reader
to Scott's excellent article~\cite{Scott83}.  Note that we always
assume that $\Sigma$ admits an orbifold cover that is a surface.  Let
$\punSig$ denote the surface obtained by removing an open neighborhood
of the orbifold points from $\Sigma$.  In most respects there is no
difference between $\Sigma$ and $\punSigma$; we will use whichever is
convenient and remark on the few subtle points as they arise.

A simple closed curve $\alpha \subset \Sigma$, avoiding the orbifold
points, is {\em inessential} if $\alpha$ bounds a disk in $\Sigma$
containing one or zero orbifold points.  The curve $\alpha$ is {\em
peripheral} if $\alpha$ is isotopic to a boundary component.  Note
that isotopies of curves are not allowed to cross orbifold points.

\begin{definition}
The {\em complex of curves} $\calC(\Sigma)$ has isotopy classes of
essential, non-peripheral curves as its vertices.  A collection of
$k+1$ distinct vertices spans a $k$--simplex if every pair of vertices
has disjoint representatives.
\end{definition}




\begin{remark}
\label{Rem:Identical}
Note that the inclusion $\punSigma \subset \Sigma$ induces an
isomorphism between $\calC(\Sigma)$ and $\calC(\punSig)$.
\end{remark}

The definition of $\calC(\Sig)$ is slightly altered when $\Sigma$ is
an annulus and also when $\punSigma$ is a once-holed torus or a
four-holed sphere.  In the last two cases the curve complex of
$\Sigma$ is the well-known {\em Farey graph}; since all curves
intersect, edges are instead placed between curves that intersect
exactly once or exactly twice, respectively.  The curve complex of an
annulus is more delicate and is defined below.

To obtain a metric, give all edges of $\calC(\Sig)$ length one and
denote distance between vertices by $d_\Sig(\cdot, \cdot)$.  It will
be enough to study only the one-skeleton of $\calC(\Sig)$, for which
we use the same notation.  This is because the one-skeleton and the
entire complex are quasi-isometric.

We begin with a simple example.

\subsection*{Puncturing}

Let $S$ be the closed surface of genus $g \geq 2$ and $\Sigma$ be the
surface of genus $g$ with one puncture. 

\begin{theorem}
\label{Thm:Puncture}
$\calC(S)$ embeds isometrically into $\calC(\Sigma)$.
\end{theorem}

As we shall see, there are uncountably many such embeddings.

\begin{proof}[Proof of \refthm{Puncture}]
Pick a hyperbolic metric on $S$.  By the Baire category theorem, the
union of geodesic representatives of simple closed curves does not
cover $S$.  (In fact, this union has Hausdorff dimension one.  See
Birman and Series~\cite{BirmanSeries85}.)
Let $\ast$ be a point in the complement and identify $\Sigma$ with
$S\setminus \{\ast\}$.  A vertex of $\calC(S)$ is then taken to its
geodesic representative, which gives an essential curve in $S
\setminus \{ \ast \}$, which is identified with a curve in $\Sigma$,
and which gives a vertex of $\calC(\Sigma)$.  This defines an
embedding $\Pi \from \calC(S) \to \calC(\Sig)$ that depends on the
choice of metric, point and identification.  Let $P \from
\calC(\Sigma) \to \calC(S)$ be the map obtained by filling the point
$\ast$. Note that $P \circ \Pi$ is the identity map.

We observe, for $a, b \in \calC(S)$ and $\alpha = \Pi(a)$, $\beta =
\Pi(b)$ that
$$
d_S(a, b) = d_\Sigma(\alpha, \beta).
$$ 
This is because $P$ and $\Pi$ send disjoint curves to disjoint curves.
Therefore, if $L \subset \calC(S)$ is a geodesic connecting $a$ and
$b$, then $\Pi(L)$ is a path in $\calC(\Sigma)$ of the same length
connecting $\alpha$ to $\beta$.  Conversely, if $\Lambda \subset
\calC(\Sigma)$ is a geodesic connecting $\alpha$ to $\beta$, then
$P(\Lambda)$ is a path in $\calC(S)$ of the same length connecting $a$
to $b$.
\end{proof}



We now turn to the main topic of the paper. 

\subsection*{Coverings}

Let $\Sigma$ and $S$ be compact connected orientable orbifolds with
negative orbifold Euler characteristic.  Let $P \from \Sigma \to S$ be
an orbifold covering map.  At a first reading it is simplest to assume
that $\Sigma$ and $S$ are both surfaces.

The covering $P$ defines a relation $\Pi \from \calC(S) \to
\calC(\Sig)$ as follows: Suppose that $b \in \calC(S)$ and $\beta \in
\calC(\Sigma)$.  Then $b$ is related to $\beta$ if and only if $\beta$
is a component of $P^{-1}(b)$, the preimage of $b$.

\begin{lemma}
The covering relation $\Pi$ is well-defined.
\end{lemma}

\begin{proof}
We will show that if $a$ is an essential non-peripheral curve then
every component of $P^{-1}(a)$ is essential and non-peripheral.  Since
$S$ has negative orbifold Euler characteristic, choose a hyperbolic
metric of finite area on $S \setminus \bdy S$.  Replace $a$ by its
geodesic representative, $a^*$.  Then $a^*$ is still a simple closed
curve, as long $a$ does not bound a disk with exactly two orbifold
points of order two.  In the latter case, $a$ collapses down to a
geodesic arc connecting the points.  Now, the lift of a geodesic
remains geodesic in the lifted metric.  The conclusion follows.
\end{proof}

We now turn to convenient piece of notation: if $A, B, \Const $ are
non-negative real numbers with $\Const > 0$ and if $A \leq \Const B +
\Const$, then we write $A \prec_\Const B$. If $A \prec_\Const B$ and
$B \prec_\Const A$, then we write $A \asymp_\Const B$.  Suppose
$\calX$ and $\calY$ are metric spaces and $f \from \calX \to \calY$ is
a relation.  We say that $f$ is a $\Const$--{\em quasi-isometric
embedding} if for all $x, x' \in \calX$ and for all $y \in f(x), y'
\in f(x')$ we have $d_\calX(x, x') \asymp_\Const d_\calY(y, y')$.  Our
goal is:

\begin{restate}{Theorem}{Thm:Cover}
The covering relation $\Pi \from \calC(S) \to \calC(\Sig)$, 
corresponding to the covering map $P \from \Sigma \to S$, is a
$\Quasi$--quasi-isometric embedding.  The constant $\Quasi$ depends
only on the topology of $S$ and the degree of the covering map $P$.
\end{restate}

\begin{remark}
Note that $\Quasi$ does not depend directly on the topology of
$\Sigma$.  When $S$ is an annulus, the degree of covering is not
determined by the topology of $\Sigma$.  Conversely, when $S$ is not
an annulus, the topology of $\Sigma$ can be bounded in terms
of the topology of $S$ and the degree of the covering.
\end{remark}

\begin{remark}
The constant $\Quasi$ may go to infinity with the degree of the
covering.  For example, any pair of distinct curves $a, b$ in a
surface $S$ may be made disjoint in some cover. In fact a cover of
degree at most $2^{d-1}$, where $d = d_S(a,b)$, will
suffice~\cite[Lemma 2.3]{Hempel01}.
\end{remark}


\begin{remark}
When $\Sigma$ is the orientation double cover of a nonorientable
surface $S$, \refthm{Cover} is due to
Masur-Schleimer~\cite{MasurSchleimer07}.
\end{remark}

The inequality $d_\Sig(\alpha, \beta) \leq d_S(a,b)$ can be easily
obtained as follows. When $\punS$ is not a once-holed torus or a
four-holed sphere, two curves in $S$ have distance one when they are
disjoint. 
But disjoint curves in $S$ have disjoint preimages in $\Sigma$.
Therefore, a path connecting $a$ to $b$ lifts to a path of equal
length connecting $\alpha$ to $\beta$ . This implies the desired
inequality in this case.  If $\punS$ is one of the special surfaces
mentioned above then two curves are at distance one when they
intersect once or twice respectively. In this case, their lifts
intersect at most $2\Degree$ times, where $\Degree$ is the degree of
the covering. Thus, the distance between the lifts of these curves is
at most $2\log_2(2\Degree) + 2$.  (See~\cite[Lemma
1.21]{Schleimer06b}.)  Therefore
$$d_\Sig(\alpha, \beta) \leq (2\log_2(2\Degree) + 2) \cdot d_S(a,b).$$

The opposite inequality is harder to obtain and occupies the bulk of
the paper.

\section{Subsurface projection}

Suppose that $\Sig$ is a compact connected orientable orbifold.  A
suborbifold $\Psi$ is {\em cleanly embedded} if every component of
$\bdy \Psi$ is either a boundary component of $\Sig$ or is an
essential non-peripheral curve in $\Sig$. All suborbifolds considered 
will be cleanly embedded.

From~\cite{MasurMinsky99}, recall the definition of the {\em
subsurface projection} relation
$$
\pi_\Psi \from \calC(\Sigma) \to \calC(\Psi),
$$ 
supposing that $\punPsi$ is not an annulus or a thrice-holed sphere.
Since $\Sigma$ has negative orbifold Euler characteristic we may
choose a complete finite volume hyperbolic metric on the interior of
$\Sigma$.  Let $\Sigma'$ be the Gromov compactification of the cover
of $\Sigma$ corresponding to the inclusion $\pi^{\orb}_1(\Psi) \to
\pi^{\orb}_1(\Sigma)$ (defined up to conjugation).  Thus $\Sigma'$ is
homeomorphic to $\Psi$; this gives a canonical identification of
$\calC(\Psi)$ with $\calC(\Sig')$. For any $\alpha \in \calC(\Sigma)$
let $\alpha'$ be the closure of the preimage of $\alpha$ in $\Sigma'$.
If every component of $\alpha'$ is properly isotopic into the boundary
then $\alpha$ is not related to any vertex of $\calC(\Psi)$; in this
case we write $\pi_\Psi(\alpha) = \emptyset$.  Otherwise, let
$\alpha''$ be a component of $\alpha'$ that is not properly isotopic
into the boundary.  Let $N$ be a closed regular neighborhood of
$\alpha'' \cup \bdy \Sigma'$.  Since $\punPsi$ is not a thrice-holed
sphere, there is a boundary component $\alpha'''$ of $N$ which is
essential and non-peripheral.  We then write $\pi_\Psi(\alpha) =
\alpha'''$.

If $\Psi$ is an annulus, then the definition of $\calC(\Psi)$ is
altered.  Vertices are proper isotopy classes of essential arcs in
$\Psi$.  Edges are placed between vertices with representatives having
disjoint interiors.  The projection map is defined as above, omitting
the final steps involving the regular neighborhood $N$.

The curve $\alpha \in \calC(\Sigma)$ {\em cuts} the suborbifold $\Psi$
if $\pi_\Psi(\alpha) \neq \emptyset$.  Otherwise, $\alpha$ {\em
misses} $\Psi$.  Suppose now that $\alpha, \beta \in \calC(\Sig)$ both
cut $\Psi$.  Define the {\em projection distance} to be
$$d_\Psi(\alpha, \beta) = d_\Psi(\pi_\Psi(\alpha), \pi_\Psi(\beta)).$$
The Bounded Geodesic Image Theorem states:

\begin{theorem}[Masur-Minsky~\cite{MasurMinsky00}]
\label{Thm:Bound}
Fix a surface $\Sigma$. There is a constant $\Bound = \Bound(\Sig)$
with the following property.  Suppose that $\alpha, \beta \in
\calC(\Sigma)$ are vertices, $\Lambda \subset \calC(\Sigma)$ is a
geodesic connecting $\alpha$ to $\beta$ and $\Omega \subsetneq \Sigma$
is a subsurface.  If $d_\Omega(\alpha, \beta) \geq \Bound$ then there
is a vertex of $\Lambda$ which misses $\Omega$. \qed
\end{theorem}

By \refrem{Identical} the Bounded Geodesic Image Theorem applies
equally well when $\Sigma$ is an orbifold.

\subsection*{Antichains}
Fix $\alpha$ and $\beta$ in $\calC(\Sigma)$ and {\em thresholds}
$\Thresh_0 > 0$ and $\Thresh_1 > 0$.  We say that a set $\calJ$ of
suborbifolds $\Omega \subsetneq \Sigma$, is a $(\Thresh_0 ,
\Thresh_1)$--{\em antichain} for $(\Sigma, \alpha, \beta)$ if $\calJ$
satisfies the following properties.
\begin{itemize}
\item
If $\Omega, \Omega' \in \calJ$ then $\Omega$ is not a strict
suborbifold of $\Omega'$.
\item
If $\Omega \in \calJ$ then $d_\Omega(\alpha, \beta) \geq \Thresh_0$. 
\item
For any $\Psi \subsetneq \Sigma$, either $\Psi$ is a suborbifold of
some element of $\calJ$ or $d_\Psi(\alpha, \beta) < \Thresh_1$.
\end{itemize}
Notice that there may be many different antichains for the given data
$(\Sigma, \alpha, \beta, \Thresh_0, \Thresh_1)$.  One particularly
nice example is when $\Thresh_0 = \Thresh_1= \Thresh$ and $\calJ$ 
is defined to be the maxima of the set 
$$\{ \Omega \subsetneq \Sigma \st 
             d_\Omega(\alpha, \beta) \geq \Thresh \}$$  
as ordered by inclusion. We call this the {\em $\Thresh$--antichain of
maxima} for $(\Sigma, \alpha, \beta)$.  By $|\calJ|$ we mean the
number of elements of $\calJ$. We may now prove:

\begin{lemma}
\label{Lem:Accumulation}
For every orbifold $\Sigma$ and for every pair of sufficiently large
thresholds $\Thresh_0, \Thresh_1$, there is an accumulation constant
$\Acc_\Sigma = \Acc(\Sigma, \Thresh_0, \Thresh_1)$ so that if $\calJ$
is an $(\Thresh_0,\Thresh_1)$--antichain for $(\Sigma, \alpha, \beta)$
then
$$d_\Sigma(\alpha, \beta) \geq |\calJ|/\Acc_\Sigma.$$
\end{lemma}

\begin{proof}
We proceed via induction: when $\calC(\punSig)$ is the Farey graph, 
$\calJ$ is the set of annuli whose core curves $\gamma$ have the property
that $d_\gamma(\alpha, \beta) \geq \Thresh_0$. In this case, 
assuming $\Thresh_0 > 3$, every such curve $\gamma$ is a vertex
of every geodesic connecting $\alpha$ to $\beta$ (see \cite[\S 4]{Minsky99}). 
Therefore the lemma holds for Farey graphs with $\Acc_\Sigma = 1$.

In the general case, let $\Error$ be a constant so that: if $\Omega
\subset \Psi \subset \Sigma$ and $\alpha', \beta'$ are the projections
of $\alpha, \beta$ to $\Psi$ then
$$|d_\Omega(\alpha, \beta) - d_\Omega(\alpha', \beta')| \leq \Error.$$ 
We take the thresholds large enough so that:
\begin{itemize}
\item
the lemma still applies to any strict suborbifold $\Psi$ with
thresholds $\Thresh_0 - \Error, \Thresh_1 + \Error$ and
\item
$\Thresh_0 \geq \Bound(\Sigma)$; thus by \refthm{Bound} for any
orbifold in $\Omega \in \calJ$ and for any geodesic $\Lambda$ in
$\calC(\Sigma)$ connecting $\alpha$ and $\beta$ there is a curve
$\gamma$ in $\Lambda$ so that $\gamma$ misses $\Omega$.
\end{itemize}
For any $\Psi \subsetneq \Sigma$ define $\Acc_\Psi = \Acc(\Psi,
\Thresh_0 - \Error, \Thresh_1 + \Error)$.

\begin{claim*}
Suppose that $\Lambda \subset \calC(\Sigma)$ is a geodesic connecting
$\alpha$ to $\beta$.  Suppose that $\gamma \in \Lambda$ is a vertex
and let $\Psi$ be a component of $\Sigma \setminus \gamma$.  The
number of elements of $\calJ_\Psi = \{ \Omega \in \calJ \st \Omega
\subsetneq \Psi \}$ is at most $\Acc_\Psi \cdot (\Thresh_1 + \Error)$.
\end{claim*}
\noindent

By the claim it will suffice to take $\Acc(\Sigma, \Thresh_0,
\Thresh_1)$ equal to
$$
2 \cdot \max \{ \Acc_\Psi \st \Psi \subsetneq \Sigma \} \cdot 
(\Thresh_1 + \Error) + 3. 
$$ 
To see this, fix a vertex $\gamma \in \Lambda$ and note that $\Sigma
\setminus \gamma$ has at most two components, say $\Psi$ and $\Psi'$.
Any element of $\calJ$ not cut by $\gamma$ is either a strict
suborbifold of $\Psi$ or $\Psi'$, an annular neighborhood of $\gamma$,
or $\Psi$ or $\Psi'$ itself.  Since every orbifold in $\calJ$ is
disjoint from some vertex of $\Lambda$, the lemma follows from the
pigeonhole principle.

All that remains is to prove the claim.  If $\Psi$ is a suborbifold of
an element of $\calJ$ then $\calJ_\Psi$ is the empty set and the claim
holds vacuously.  Thus we may assume that
$$
d_\Psi(\alpha, \beta) < \Thresh_1.
$$ 
Let $\alpha'$ and $\beta'$ be the projections of $\alpha$ and $\beta$ 
to $\Psi$.  From the definition of $\Error$, $\calJ_\Psi$ is a 
$(\Thresh_0 - \Error, \Thresh_1 + \Error)$--antichain for 
$\Psi, \alpha'$ and $\beta'$. Thus,
$$\Thresh_1 > d_\Psi(\alpha, \beta) \geq d_\Psi(\alpha', \beta') -
\Error \geq |\calJ_\Psi|/\Acc_\Psi - \Error,$$ 
with the last inequality following by induction. Hence, 
\begin{equation*}
\Thresh_1 + \Error \geq |\calJ_\Psi|/\Acc_\Psi. \qedhere
\end{equation*}
\end{proof}

\section{\Teich space}

For this section, we take $\Sigma$ to be a surface.  Let $\calT(\Sig)$
denote the \Teich space of $\Sigma$: the space of complete hyperbolic
metrics on the interior of $\Sigma$, up to isotopy.  For background,
see \cite{Bers60, Gardiner87}.  


There is a uniform upper bound on the length of the shortest closed
curve in any hyperbolic metric on $\Sigma$. For any metric $\sigma$ on
$\Sigma$, a curve $\gamma$ has \emph{bounded length} in $\sigma$ if
the length of $\gamma$ in $\sigma$ is less than this constant.  Let
$\ep_0>0$ be a constant such that, for curves $\gamma$ and $\delta$,
if $\gamma$ has bounded length in $\sigma$ and $\delta$ has a length
less than $\ep_0$ then $\gamma$ and $\delta$ have intersection number
zero.

Suppose that $\alpha$ and $\beta$ are vertices of $\calC(\Sig)$.  Fix
metrics $\sigma$ and $\tau$ in $\calT(\Sig)$ so that $\alpha$ and
$\beta$ have bounded length at $\sigma$ and $\tau$ respectively.  Let
$\Gamma \from [t_\sigma, t_\tau] \to \calT(S)$ be a geodesic
connecting $\sigma$ to $\tau$.  For any curve $\gamma$ let
$l_t(\gamma)$ be the length of its geodesic representative in the
hyperbolic metric $\Gamma(t)$.  
The following theorems are consequences of Theorem~6.2 and Lemma~7.3
in \cite{Rafi05}.

\begin{theorem}[\cite{Rafi05}]
\label{Thm:Large->Short}
For $\ep_0$ as above there exists a threshold $\Thresh_{\rm min}$ such
that, for a strict subsurface $\Omega$ of $\Sigma$, if
$d_\Omega(\alpha, \beta) \geq \Thresh_{\rm min}$ then there is a time
$t_\Omega$ so that the length of each boundary component of $\Omega$
in $\Gamma(t_\Omega)$ is less than $\ep_0$. \qed
\end{theorem}

\begin{theorem}[\cite{Rafi05}] 
\label{Thm:Short->Large}
For every threshold $\Thresh_1$, there is a constant $\ep_1$ such that
if $l_t(\gamma) \leq \ep_1$, for some curve $\gamma$ and for some time
$t$, then there exists a subsurface $\Psi$ disjoint from $\gamma$ such
that $d_\Psi(\alpha, \beta) \geq \Thresh_1$. \qed
\end{theorem}



The {\em shadow} of the \Teich geodesic $\Gamma$ inside of
$\calC(\Sigma)$ is the set of curves $\gamma$ so that $\gamma$ has
bounded length in $\Gamma(t)$ for some $t\in[t_\sigma, t_\tau]$. The
following is a consequence of the fact that the shadow is an
unparameterized quasi-geodesic.  (See Theorem~2.6 and then apply
Theorem~2.3 in \cite{MasurMinsky99}.)

\begin{theorem}[\cite{MasurMinsky99}]
\label{Thm:BackTrack}
The shadow of a \Teich geodesic inside of $\calC(\Sigma)$ does not
backtrack and so satisfies the reverse triangle inequality.  That is,
there exists a backtracking constant $\BackTrack = \BackTrack(\Sigma)$
such that if $t_\sigma \leq t_0 \leq t_1 \leq t_2 \leq t_\tau$ and if
$\gamma_i$ has bounded length in $\Gamma(t_i)$, $i = 0, 1, 2$ then
\begin{equation*} 
d_\Sigma(\gamma_0, \gamma_2) \geq 
           d_\Sigma(\gamma_0, \gamma_1) + 
           d_\Sigma(\gamma_1, \gamma_2) - \BackTrack. \qedhere 
\end{equation*}
\end{theorem}

We say that $\Gamma(t)$ is $\ep$--{\em thick} if the shortest closed
geodesic $\gamma$ in $\Gamma(t)$ has a length of at least $\ep$. 

\begin{lemma}
\label{Lem:Thick}
For every $\ep > 0$ there is a progress constant $\Progress > 0$ so that
if $t_\sigma \leq t_0 \leq t_1 \leq t_\tau$, if $\Gamma(t)$ is
$\ep$--thick at every time $t \in [t_0, t_1]$, and if $\gamma_i$ has
bounded length in $\Gamma(t_i)$ ($i = 0, 1$) then
\[
d_\Sigma(\gamma_0, \gamma_1) \asymp_\Progress t_1 - t_0.
\qedhere \]
\end{lemma}

\begin{proof}
As above, using Theorem~6.2 and Lemma~7.3 in \cite{Rafi05} and the
fact that $\Gamma(t)$ is $\ep$--thick at every time $t \in [t_0,
t_1]$, we can conclude that $d_\Omega(\gamma_0, \gamma_1)$ is
uniformly bounded for any strict subsurface of $\Omega$ of $\Sigma$.
The lemma is then a consequence of Theorem~1.1 and Remark~5.5 in
\cite{Rafi06}.  (Referring to the statement and notation
of~\cite[Theorem~1.1]{Rafi06}: Extend $\gamma_i$ to a short marking
$\mu_i$.  Take $k$ large enough such that the only non-zero term in
the right hand side of \cite[Equation (1)]{Rafi06} is $d_\Sigma(\mu_0,
\mu_1)$.)
\end{proof}

In general the geodesic $\Gamma$ may stray into the thin part of
$\calT(S)$.  We take $\Gamma^{\geq \ep}$ to be the set of times in the
domain of $\Gamma$ which are $\ep$--thick.  Notice that $\Gamma^{\geq
\ep}$ is a union of closed intervals.  Let $\LongThick$ be the union
of intervals of $\Gamma^{\geq \ep}$ which have length at least
$\Length$.  We use $|\LongThick|$ to denote the sum of the lengths of
the components of $\LongThick$.

\begin{lemma} 
\label{Lem:LongThick}
For every $\ep$ there exists $\Length_0$ such that if $\Length \geq
\Length_0$, then
$$d_\Sigma(\alpha, \beta) \geq |\LongThick|/2\Progress.$$ 
\end{lemma}

\begin{proof}
Pick $\Length_0$ large enough so that, for $\Length \geq \Length_0$,
$$
(\Length/2\Progress) \geq \Progress + 2\BackTrack.
$$

Let $\LongThick$ be the union of intervals $[t_i, s_i]$, $i = 1,
\ldots, m$.  Let $\gamma_i$ be a curve of bounded length in
$\Gamma(t_i)$ and $\delta_i$ be a curve of bounded length in
$\Gamma(s_i)$.

By \refthm{BackTrack} we have
$$d_\Sig(\alpha, \beta) \geq 
     \left( \sum_i d_\Sig(\gamma_i, \delta_i) \right) - 2m
     \BackTrack.$$
From \reflem{Thick} we deduce
$$d_\Sig(\alpha, \beta) \geq 
     \left( \sum_i \frac{1}{\Progress} (s_i - t_i) 
     - \Progress \right) - 2m \BackTrack.$$
Rearranging, we find
$$d_\Sig(\alpha, \beta) \geq 
     \frac{1}{\Progress} |\LongThick| - m(\Progress + 2\BackTrack).$$
Thus, as desired:
\[ d_\Sig(\alpha, \beta) \geq 
     \frac{1}{2\Progress} |\LongThick|. \qedhere \]
\end{proof}

\section{An estimate of distance}

Again, take $\Sigma$ to be a surface.  In this section we provide the
main estimate for $d_\Sigma(\alpha, \beta)$.  Let $\ep_0$ be as
before. We choose thresholds $\Thresh_0 \geq \Thresh_{\rm min}$ (see
Theorem~\ref{Thm:Large->Short}) and $\Thresh_1$ so that
\reflem{Accumulation} holds. Let $\ep_1$ be the constant provided in
\reflem{Thick} and let $\ep>0$ be any constant smaller than $\min \{
\ep_0, \ep_1 \}$. Finally, we pick $\Length_0$ such that
\reflem{LongThick} holds and that $\Length_0/2\Progress > 4$.  Let
$\Length$ be any length larger than $\Length_0$.

\begin{theorem}
\label{Thm:Estimate}
Let $\Thresh_0$, $\Thresh_1$, $\ep$ and $\Length$ be constants chosen
as above.  There is a constant $\Comp = \Comp(\Sigma, \Thresh_0,
\Thresh_1, \ep , \Length)$ such that for any curves $\alpha$ and
$\beta$, any $(\Thresh_0, \Thresh_1)$--antichain $\calJ$ and any
\Teich geodesic $\Gamma$, chosen as above, we have:
$$d_\Sigma(\alpha, \beta) \asymp_\Comp |\calJ| + |\LongThick|.$$ 
\end{theorem}

\begin{proof}
For $\Comp \geq 2 \cdot \max(\Acc, 2\Progress)$, the inequality
$$d_\Sigma(\alpha, \beta) \succ_\Comp |\calJ| + |\LongThick|$$ follows
from Lemmas~\ref{Lem:Accumulation} and~\ref{Lem:LongThick}.  It
remains to show that
$$d_\Sig(\alpha, \beta) \prec_\Comp |\calJ|+ |\LongThick|.$$

For each $\Omega \in \calJ$ fix a time $t_\Omega \in [t_\sigma,
t_\tau]$ so that all boundary components of $\Omega$ are
$\ep_0$--short in $\Gamma(t_\Omega)$ (see
Theorem~\ref{Thm:Large->Short}). Let $\calE$ be the union:
$$
\Big\{t_\Omega \ \Big|\ \Omega \in \calJ, \ t_\Omega \not \in \LongThick \Big\}
\cup 
\Big\{\bdy I \ \Big|\ \text{$I$ a component of  $\LongThick$} \Big\}.
$$
We write $\calE = \{ t_0,
\ldots, t_n \}$, indexed so that $t_i < t_{i+1}$.

\begin{claim*}
The number of intervals in $\LongThick$ is at most $|\calJ| + 1$.
Hence, $|\calE| \leq 3 |\calJ| + 1$. 
\end{claim*}

\begin{proof}
There is at least one moment between any two consecutive intervals $I,
J \subset \LongThick$ when some curve $\gamma$ becomes $\ep$--short
(and hence $\ep_1$--short).  Therefore, by
Theorem~\ref{Thm:Short->Large}, $\gamma$ is disjoint from a subsurface
$\Psi$ where $d_\Psi(\alpha, \beta)\geq \Thresh_1$. Since $\calJ$ is
an $(\Thresh_0, \Thresh_1)$--antichain, $\Psi$ is a subsurface of some
element $\Omega \in \calJ$.  It follows that $d_\Sigma(\gamma, \bdy
\Omega) \leq 2$. This defines a one-to-many relation from pairs of
consecutive intervals to $\calJ$.  To see the injectivity consider
another such pair of consecutive intervals $I'$ and $J'$ and a
corresponding curve $\gamma'$ and subsurface $\Omega'$. By
\reflem{Thick}, $d_\Sigma(\gamma, \gamma') \geq \Length/2\Progress >
4$ and therefore $\Omega$ is not equal to $\Omega'$.
\end{proof}

Let $\gamma_i$ be a curve of bounded length in $\Gamma(t_i)$.
\begin{claim*}
$$d_\Sig(\gamma_i, \gamma_{i+1}) \leq
\begin{cases}
\Progress(t_{i+1} - t_i) + \Progress,  
    &\quad \text{if $[t_i, t_{i+1}] \subset \LongThick$} \\
2\BackTrack + \Progress\Length + \Progress + 2,
    &\quad \text{otherwise.}
\end{cases}$$
\end{claim*}

\begin{proof}
The first case follows from \reflem{Thick}.  So suppose that the
interior of $[t_i, t_{i+1}]$ is disjoint from the interior of
$\LongThick$.  

We define sets $I_+, I_- \subset [t_i, t_{i+1}]$ as
follows: A point $t \in [t_i, t_{i+1}]$ lies in $I_-$ if
\begin{itemize}
\item
there is a curve $\gamma$ which is $\ep$--short in $\Gamma(t)$ and
\item
for some $\Omega \in \calJ$, so that $d_\Sig(\bdy \Omega, \gamma) \leq
2$, we have $t_\Omega \leq t_i$.
\end{itemize}
If instead $t_\Omega \geq t_{i+1}$ then we
place $t$ in $I_+$.  Finally, we place $t_i$ in $I_-$ and $t_{i+1}$ in
$I_+$. 

Notice that if $\Omega \in \calJ$ then $t_\Omega$ does not lie in the
open interval $(t_i, t_{i+1})$.  It follows that every $\ep$--thin
point of $[t_i, t_{i+1}]$ lies in $I_-$, $I_+$, or both.  If $t\in
I_-$ and $\gamma$ is the corresponding $\ep$--short curve then
$d_\Sigma(\gamma_i, \gamma)\leq \BackTrack + 2$.  This is because
either $t = t_i$ and so $\gamma$ and $\gamma_i$ are in fact disjoint,
or there is a surface $\Omega \in \calJ$ as above with
$$ 
2 \geq d_\Sigma(\bdy \Omega, \gamma) \geq 
   d_\Sigma(\gamma_i, \gamma) - \BackTrack,
$$ 
Similarly if $t \in I_+$ then $d_\Sigma(\gamma_{i+1}, \gamma)\leq
\BackTrack + 2$.

If $I_+$ and $I_-$ have non-empty intersection then
$d_\Sigma(\gamma_i, \gamma_{i+1}) \leq 2B + 4$ by the triangle
inequality.

Otherwise, there is an interval $[s,s']$ that is $\ep$--thick, has
length less than $\Length$ such that $s \in I_-$ and $s' \in I_+$.
Let $\gamma$ and $\gamma'$ be the corresponding short curves in
$\Gamma(s)$ and $\Gamma(s')$. Thus
$$
d_\Sigma(\gamma_i, \gamma) \leq \BackTrack + 2
\quad\text{and}\quad d_\Sigma(\gamma', \gamma_{i+1}) \leq \BackTrack +2.
$$
We also know from \reflem{Thick} that
$$d_\Sigma(\gamma, \gamma') \leq \Progress\Length+\Progress.$$
This finishes the proof of our claim. 
\end{proof}

It follows that
\begin{align*}
d_\Sig(\alpha, \beta) 
   & \leq d_\Sig(\gamma_0, \gamma_{1}) + \ldots + 
                        d_\Sig(\gamma_{n-1}, \gamma_{n}) \\
   & \leq |\calE| (2\BackTrack + \Progress \Length + \Progress +2 ) 
                       + \Progress |\LongThick| + |\calE| \Progress \\
                                 & \prec_\Comp |\calJ| +  |\LongThick|,
\end{align*}
for an appropriate choice of $\Comp$.  This proves the theorem. 
\end{proof}

\section{Symmetric curves and surfaces}

Let $P \from \Sig \to S$ be an orbifold covering map.

\begin{definition}
A curve $\alpha \in \calC(\Sig)$ is {\em symmetric} if there is a
curve $a \in \calC(S)$ so that $P(\alpha) = a$.  We make the same
definition for a suborbifold $\Omega \subset \Sigma$ lying over a
suborbifold $Z \subset S$.
\end{definition}

For the rest of this section, fix symmetric curves $\alpha$ and
$\beta$.  Pick $x,y \in \calT(\punS)$ so that $a = P(\alpha)$ has
bounded length in $x$ and $b = P(\beta)$ is bounded in $y$.  Let $G
\from [t_x, t_y] \to \calT(\punS)$ be the \Teich geodesic connecting
$x$ to $y$.  For every $t \in [t_x, t_y]$ let $q_t$ be the terminal
quadratic differential of the \Teich map from $G(t_x)$ to $G(t)$.  We
lift $q_t$ to the surface $P^{-1}(\punS)$, fill the punctures not
corresponding to orbifold points and so obtain a parameterized family
$\theta_t$ of quadratic differentials on $\punSig$.  Notice that
$\theta_t$ is indeed a quadratic differential: suppose that $p \in S$
is a orbifold point and $q_t$ has a once-pronged singularity at $p$.
For every regular point $\pi$ in the preimage of $p$ the differential
$\theta_t$ has at least a twice-pronged singularity at $\pi$.

Uniformize the associated flat structures to obtain hyperbolic metrics
on $\punSig$.  This gives a path $\Gamma \from [t_x, t_y] \to
\calT(\punSig)$.  The path $\Gamma$ is a geodesic in $\calT(\punSigma)$.
This is because, for $t, s \in [t_x, t_y]$, the \Teich map from $G(t)$
to $G(s)$ has Beltrami coefficient ${\sf k} \, |q|/q$ where $q$ is an
integrable holomorphic quadratic differential in $G(t)$. This map
lifts to a map from $\Gamma(t)$ to $\Gamma(s)$ with Beltrami
coefficient ${\sf k} \, |\theta|/\theta$, where the quadratic
differential $\theta$ is the pullback of $q$ to $\Gamma(t)$. That is,
the lift of the \Teich map from from $G(t)$ to $G(s)$ is the \Teich
map from $\Gamma(t)$ to $\Gamma(s)$ with the same quasi-conformal
constant. Therefore, as is well-known, the distance in $\calT(\punS)$
between $G(t)$ and $G(s)$ equals the distance in $\calT(\punSigma)$
between $\Gamma(t)$ and $\Gamma(s)$.

\begin{proposition}[Proposition 3.7 \cite{Rafi06}]
For any $\ep$, there is a constant $\NoProgress$ such that the
following holds.  Assume that, for all $t \in [r,s]$, there is a
component of $\bdy \Omega$ whose length in $\Gamma(t)$ is larger than
$\ep$. Suppose $\gamma$ has bounded length in $\Gamma(r)$ and $\delta$
has bounded length in $\Gamma(s)$. Then
\[ 
d_\Omega(\gamma, \delta) \leq \NoProgress. \qedhere 
\]
\end{proposition}

\begin{lemma} 
\label{Lem:Symmetric}
For $\ep$ small enough, $\NoProgress$ as above and any suborbifold
$\Omega \subset \Sig$, if $d_\Omega(\alpha, \beta) \geq 2\NoProgress +
1$, then $\Omega$ is symmetric.
\end{lemma}

\begin{proof}
Consider the first time $t^-$ and last time $t^+$ that the boundary of
$\Omega$ is $\ep$--short. Since every component of $\bdy \Omega$ is
short in $\Gamma(t^\pm)$, so is the image $P(\bdy \Omega)$ in
$G(t^\pm)$.  Therefore, all components of the image are simple.  (This
is a version of the Collar Lemma.  For example, see
\cite[Theorem~4.2.2]{Buser92}.)  It follows that the boundary of
$\Omega$ is symmetric. This is because choosing $\ep$ small enough
will ensure that curves in $P^{-1} ( P( \Omega))$ have bounded length
at both $t^-$ and $t^+$. (The length of each is at most the degree of
the covering map times $\ep$.)  If any such curve $\gamma$ intersects
$\Omega$ we have $d_\Omega(\gamma, \alpha) \leq \NoProgress$ and
$d_\Omega(\gamma, \beta) \leq \NoProgress$, contradicting the
assumption $d_\Omega(\alpha, \beta) \geq 2\NoProgress + 1$.  Thus, the
suborbifold $\Omega$ is symmetric.
\end{proof}

\section{The quasi-isometric embedding}

We now  prove the main theorem:

\begin{theorem}
\label{Thm:Cover}
The covering relation $\Pi \from \calC(S) \to \calC(\Sig)$, 
corresponding to the covering map $P \from \Sigma \to S$, is a
$\Quasi$--quasi-isometric embedding.  The constant $\Quasi$ depends
only on the topology of $S$ and the degree of the covering map $P$.
\end{theorem}

\begin{proof}
As mentioned before, we only need to show that
$$
d_\Sig(\alpha, \beta) \succ_\Quasi d_S(a,b).
$$ 
Suppose that $\Degree$ is the degree of the covering. We prove the
theorem by induction on the complexity of $S$. In the case where $S$
is an annulus without orbifold points, the cover $\Sigma$ is also an
annulus and the distances in $\calC(\Sigma)$ and $\calC(S)$ are equal
to the intersection number plus one. But, in this case,
$$
i(\alpha, \beta) \geq i(a,b) /\Degree.
$$ 
Therefore, the theorem is true with $\Quasi = \Degree$.

Now assume the theorem is true for all strict suborbifolds of $S$.
Let $\Quasi'$ be the largest constant of quasi-isometry necessary for
such suborbifolds.  Choose the threshold $\Thresh$, constant $\ep$ and
length $\Length$ such that Theorem~\ref{Thm:Estimate} holds for both
the data $(S, \Thresh, \Thresh, \ep, \Length)$ as well as $(\Sigma,
(\Thresh/\Quasi') - \Quasi', \Thresh, \ep, \Length)$. We also assume
that $\Thresh \geq 2 \NoProgress + 1$.  All of the constants depend
only on the topology of $S$ and the degree $\Degree$, because these
last two bound the topology of $\Sigma$.

Let $\calJ_S$ be the $\Thresh$--antichain of maxima for $S, a$ and $b$
and let $\calJ_\Sigma$ be the set of components of preimages of
elements of $\calJ_S$.
\begin{claim*}
The set $\calJ_\Sigma$ is a $((\Thresh/\Quasi') - \Quasi',
\Thresh)$--antichain for $(\Sigma, \alpha, \beta)$.
\end{claim*}
We check the conditions for being an antichain. Since elements of
$\calJ_S$ are not subsets of each other, the same holds for their
preimages.  The condition $d_\Omega(\alpha, \beta) \geq
(\Thresh/\Quasi') - \Quasi'$ is the induction hypothesis. Now suppose
$\Psi \subset \Sigma$ with $d_\Psi(\alpha, \beta) \geq \Thresh$.  By
\reflem{Symmetric}, $\Psi$ is symmetric.  That is, it is a component
of the preimage of an orbifold $Y\subset S$ and
$$
d_Y(a, b) \geq d_\Psi(\alpha, \beta) \geq \Thresh.
$$ 
This implies that $Y \subset Z$ for some $Z \in \calJ_S$. Therefore,
taking $\Omega$ to be the preimage of $Z$, we have $\Psi \subset
\Omega \in \calJ_\Sigma$.  This proves the claim.

Hence, there are constants $\Comp$ and $\Comp'$ such that 
\begin{align*}
d_S(a,b) &\asymp_\Comp |\calJ_S| + |G(\ep, \Length)|,\\
\intertext{and}
d_\Sigma(\alpha,\beta) 
&\asymp_{\Comp'} |\calJ_\Sigma| + |\Gamma(\ep, \Length)|.
\end{align*}
Note that $|\calJ_S| \leq \Degree |\calJ_\Sigma|$ as a suborbifold of
$S$ has at most $\Degree$ preimages.  Note also that $|G(\ep,
\Length)| \leq |\Gamma(\ep, \Length)|$ because $\Gamma(t)$ is at least
as thick as $G(t)$.  Therefore
$$
d_S(a,b) \prec_\Quasi d_\Sigma(\alpha, \beta),
$$
for $\Quasi = \Degree \, \Comp \, \Comp'$. This finishes the proof.
\end{proof}

\section{An application to mapping class groups}

Suppose that $P \from \Sigma \to S$ is an orbifold covering map, and
$\chi^{\orb}(S) < 0$.  Let $\MCG(\Sig)$ be the orbifold mapping class
group of $\Sig$: isotopy classes of homeomorphisms of $\Sigma$
restricting to the identity on $\bdy S$ and respecting the set of
orbifold points and their orders.  Here all isotopies must fix all
boundary components and all orbifold points.  As an application of
\refthm{Cover} we prove the following theorem:

\begin{theorem}
\label{Thm:MCG}
The covering $P$ induces a quasi-isometric embedding 
$$
\Pi_* \from \MCG(S) \to \MCG(\Sigma).
$$
\end{theorem}

\begin{remark}
Assume that $P \from \Sigma \to S$ is a regular orbifold cover with
finite deck group $\Delta < \MCG(\Sigma)$.  Let $N(\Delta)$ be the
normalizer of $\Delta$ inside of $\MCG(\Sigma)$ and let $M < \MCG(S)$
be the finite index subgroup of mapping classes that lift.  MacLachlan
and Harvey ~\cite[Theorem~10]{MacLachlanHarvey75}: showed that there
is a short exact sequence:
\begin{equation*}
1 \to \Delta \to N(\Delta) \to M \to 1.
\end{equation*}
It is not difficult to show that any set-theoretic section of the
MacLachlan-Harvey map is bounded distance from the quasi-isometric
embedding constructed by the proof of \refthm{MCG}.
\end{remark}

\begin{remark}
\refthm{MCG} can also be compared to the following statement: 
Suppose that a subsurface $Z \subset S$ is cleanly embedded and 
$S \setminus Z$ has no annuli components.  Then the
inclusion $Z \to S$ induces a quasi-isometric
embedding on mapping class groups.  This follows directly from the
summation formula of Masur and Minsky (See~\cite{MasurMinsky00},
Theorems 7.1, 6.10, and 6.12) and was independently obtained by 
Hamenst\"adt \cite[Theorem~B, Corollary~4.6]{Hamenstaedt05a}.
\end{remark}

\begin{proof}[Proof of \refthm{MCG}]
Choose, for the remainder of the proof, a {\em marking} of $S$: 
a collection, $m$, of curves which fill $S$. 
Let $\mu = \Pi(m)$ be the lift of $m$ to $\Sigma$. Since disks and punctured
disks lift to disks and punctured disks the curves of $\mu$
fill $\Sigma$. Hence $\mu$ is a marking on $\Sigma$. 

We construct $ \Pi_\star$ as follows: Let $x$ be an element of
$\MCG(S)$. Define $\xi=\Pi_\star(x)$ to be any element of
$\MCG(\Sigma)$ so that $\xi(\mu)$ intersects $\Pi(x(m))$ an {\em a
priori} bounded number of times. Such a map $\xi$ always exists and,
for a given bound on intersection, there are only finitely many
possibilities for such a map (see \cite{MasurMinsky00}).

Let $T$ be a generating set for $\MCG(S)$ and $\Theta$ be a generating
set for $\MCG(\Sigma)$. Let $||x||_T$ and $||\xi||_\Theta$ denote the
word lengths of $x$ and $\xi$ with respect to $T$ and $\Theta$
respectively.  To prove the proposition it is sufficient to show that,
for $\xi = \Pi_\star(x)$,
\begin{equation} \label{Eqn:WordLength}
||x||_T \asymp_{\sf W} ||\xi||_\Theta, 
\end{equation}
where ${\sf W}$ is a constant that does not depend on $x$. 

By \cite{MasurMinsky00}, Theorems 7.1, 6.10, and 6.12, we have
\begin{equation} \label{Eqn:LengthInS}
||x||_T \asymp_{\sf W_1} \sum \Big[d_Z\big(m, x(m)\big) \Big]_{\sf k_1}.
\end{equation}
Here the sum ranges over all sub-orbifolds $Z \subset S$.  The
constant ${\sf W_1}$ depends on ${\sf k_1}$ which in turn depends on
our choice of the marking $m$ and the generating set $T$.  However,
all of the choices are independent of the group element $x$.  Finally,
$[r]_{\sf k} = r$ if $r \geq {\sf k}$ and $[r]_{\sf k} = 0$ if $r <
{\sf k}$.

As above, after fixing a large enough constant ${\sf k_2}$ (see below) 
and an appropriate ${\sf W_2}$, we have
\begin{equation}
||\xi||_\Theta \asymp_{\sf W_2}  \sum
\Big[d_\Omega\big(\mu, \xi(\mu) \big) \Big]_{\sf k_2}.
\end{equation}
But $\xi(\mu)$ and $\Pi(x(m))$ have bounded intersection.  Therefore,
their projection distance in every subsurface $\Omega$ is {\em a
priori} bounded. Hence we can write
\begin{equation} \label{Eqn:LengthInSig}
||\xi||_\Theta \asymp_{\sf W_3}  \sum
\Big[d_\Omega\big(\mu, \Pi(x(m)) \big) \Big]_{\sf k_2},
\end{equation}
for a slightly larger constant ${\sf W}_3$.

We prove equation \eqref{Eqn:WordLength} 
by comparing the terms of the the right hand side of \eqref{Eqn:LengthInS} 
with those on the right hand side of \eqref{Eqn:LengthInSig}.
Note that $\mu=\Pi(m)$ is a union of symmetric orbits and the same holds for 
$\Pi(x(m))$.  Therefore, we can choose ${\sf k_2}$ large enough such that if 
$d_\Omega(\mu, g(\mu))$ is larger than ${\sf k_2}$ 
then $\Omega$ is itself symmetric (see \reflem{Symmetric}).  
Taking $Z = P(\Omega)$, it follows from \refthm{Cover} that
$$
d_Z(m, x(m)) \asymp d_\Omega(\mu, \Pi(x(\mu))).
$$ 
On the other hand, \refthm{Cover} also tells us that large projection
distance in any $Z \subset S$ implies large projection distance in all
the components of the pre-image of $Z$.  Therefore, there is a
finite-to-one correspondence between the surfaces that appear in
\eqref{Eqn:LengthInSig} and in \eqref{Eqn:LengthInS} and the
corresponding projection distances are comparable. We conclude that
$||x||_T \asymp_{\sf W} ||\xi||_\Theta$ for some ${\sf W}$. This
finishes the proof.
\end{proof}

As a further application, let $\Sigma$ be the closed orientable
surface of genus $g$ and let $\phi \from \Sigma \to \Sigma$ be a
hyperelliptic involution.  Let $S = \Sigma/\phi$ and let $P \from \Sigma
\to S$ be the induced orbifold cover.  Birman and
Hilden~\cite{BirmanHilden73} provide a short exact sequence:
$$
1 \to \langle \phi \rangle \to N(\phi) \to \MCG(S) \to 1.
$$ 
Notice that $\MCG(S)$ is the spherical braid group on $2g+2$ strands.
As in \refthm{MCG} we immediately have:

\begin{corollary}
\label{Cor:Luis}
A section of the Birman-Hilden map induces a quasi-isometric embedding
of the spherical braid group on $2g+2$ strands into the mapping class
group of the closed surface of genus $g$. \qed
\end{corollary}

\bibliographystyle{alpha} 
\bibliography{bibfile}
\end{document}